\documentclass[12pt]{article}

\usepackage{epsfig}
\usepackage{amssymb}
\usepackage{amsmath}

\def\shadowbox{\hbox{\rule[-0.0ex]{0.1ex}{1.2ex}%
\hspace{-0.1ex}\rule[-0.0ex]{1.2ex}{0.1ex}%
\hspace{0.0ex}\rule[-0.0ex]{0.1ex}{1.2ex}\hspace{-1.3ex}%
\rule[1.15ex]{1.25ex}{0.1ex}\hspace{-0.0ex}\rule[-0.25ex]{0.3ex}{1.1ex}%
\hspace{-1.2ex}\rule[-0.25ex]{1.1ex}{0.25ex}}}
\def\qed{\ifmmode \hbox{\hfill\shadowbox}
     \else \hphantom{x}\hfill\shadowbox \fi}
\def\QED{\mbox{\phantom{m}}\nolinebreak\hfill$\,\shadowbox$}
\def\proof{\noindent{\bf Proof: }}

\newtheorem{theorem}{Theorem}[section]
\newtheorem{lemma}[theorem]{Lemma}
\newtheorem{definition}[theorem]{Definition}
\newtheorem{proposition}[theorem]{Proposition}
\newtheorem{corollary}[theorem]{Corollary}

\def\Cst{{\mathbb C}}

\def\Qst{{\mathbb Q}}
\def\Rst{{\mathbb R}}
\def\Zst{{\mathbb Z}}
\def\Lsp{{\boldsymbol L}}
\def\Ban{{\cal B}}
\def\Win{{\cal W}}
\def\eps{\varepsilon}
\def\phi{\varphi}
\def\Cst{{\mathbb C}}

\def\Qst{{\mathbb Q}}
\def\Rst{{\mathbb R}}

\def\Zst{{\mathbb Z}}

\def\Lsp{{\boldsymbol L}}
\def\lsp{{\boldsymbol\ell}}
\def\Ltsp{{\Lsp^2}}
\def\LtR{{\Lsp^2(\Rst)}}
\def\ltZ{{\lsp^2(\Zst)}}
\def\ltsp{{\lsp^2}}
\def\Ssp{{\boldsymbol S}}
\def\SOsp{{\Ssp_0}}

\def\adj{\operatorname{adj}}
\def\argmin{\operatorname{arg\, min}}
\def\trace{\operatorname{trace}}

\def\Re{\operatorname{Re}}

\def\aet{{\text{a.e.}\,t}}
\def\aen{{\text{a.e.}\, \nu}}

\def\aetn{{\text{a.e.}\,t, \nu}}

\newcommand{\Ti}{T^{-\frac{1}{2}}}
\newcommand{\TKi}{T^{(K)}}
\newcommand{\Tni}{T_N^{-\frac{1}{2}}}
\newcommand{\TKni}{T_N^{(K)}}
\newcommand{\Phih}{\Phi^h (t,\nu)}
\newcommand{\Phig}{\Phi^g (t,\nu)}
\newcommand{\Phiha}{\big(\Phi^h(t,\nu)\big)^{\ast}}
\newcommand{\Phiga}{\big(\Phi^g(t,\nu)\big)^{\ast}}

\newcommand{\SI}{S^{-1}}
\newcommand{\SQI}{S^{-\frac{1}{2}}}

\newcommand{\gkl}{g_{k/b,l/a}}

\newcommand{\gklp}{g_{k'/b,l'/a}}
\newcommand{\hkl}{h_{k/b,l/a}}
\newcommand{\ho}{{h^{0}}}
\newcommand{\gd}{\gamma^{0}}
\newcommand{\gdnm}{\gd_{na,mb}}

\newcommand{\gnm}{g_{na,mb}}
\newcommand{\gsnm}{g_{n,m}}
\newcommand{\hnm}{h_{na,mb}}

\newcommand{\honm}{\ho_{\!\!na,mb}}
\newcommand{\hosnm}{\ho_{\!\!n,m}}
\newcommand{\homm}{\ho_{\!\!m}}

\newcommand{\kb}{\frac{k}{b}}
\newcommand{\Mg}{M_g}
\newcommand{\Mga}{M^{\ast}_g}
\newcommand{\Hg}{H_g}
\newcommand{\Hga}{H^{\ast}_g}
\newcommand{\Mho}{M_{\ho}}
\newcommand{\Mgka}{M^{\ast}_{g,K}}
\newcommand{\Mh}{M_{h}}
\newcommand{\Mha}{M_{h}^{\ast}}
\newcommand{\Mgk}{M_{g,K}}
\newcommand{\Mhk}{M_{h,K}}

\newcommand{\Mk}{M_{K}}
\newcommand{\Mka}{M_{K}^{\ast}}
\newcommand{\Ug}{U_g}
\newcommand{\Uga}{U^{\ast}_g}
\newcommand{\Ugkl}{U_{g,K,L}}
\newcommand{\Ugkla}{U^{\ast}_{g,K,L}}
\newcommand{\Uh}{U_h}
\newcommand{\Uha}{U^{\ast}_h}
\newcommand{\Uhkl}{U_{h,K,L}}
\newcommand{\Uhkla}{U^{\ast}_{h,K,L}}

\def\toinf{{\rightarrow \infty }}

\newcommand{\gabframe}{(\gnm)_{n,m \in \Zst}}
\newcommand{\tightframe}{(\honm)_{n,m \in \Zst}}

\begin{document}

\title{Characterization and computation of canonical tight windows
for Gabor frames}

\author{A.J.E.M.~Janssen\thanks{Philips Research Laboratories WY-81, 5656 AA 
Eindhoven, The Netherlands; Email: a.j.e.m.janssen@philips.com} and Thomas 
Strohmer\thanks{Department of Mathematics, University of California, Davis, 
CA 95616-8633, USA; Email: strohmer@math.ucdavis.edu. T.S. acknowledges 
support from NSF grant 9973373.}}

\date{}
\maketitle



\begin{abstract}
Let $(\gnm)_{n,m\in\Zst}$ be a Gabor frame for $\LtR$ for given window $g$.
We show that the window $\ho=\SQI g$ that 
generates the canonically associated tight Gabor frame minimizes
$\|g-h\|$ among all windows $h$ generating a normalized tight Gabor frame.
We present and prove versions of this result in the time domain, the
frequency domain, the time-frequency domain, and the Zak transform domain,
where in each domain the canonical $\ho$ is expressed using functional 
calculus for Gabor frame operators. Furthermore, we derive a Wiener-Levy 
type theorem for rationally oversampled Gabor frames.
Finally, a Newton-type method for a fast numerical calculation of $\ho$ 
is presented. We analyze the convergence behavior of this method and 
demonstrate the efficiency of the proposed algorithm by some numerical 
examples.
\end{abstract}



\section{Introduction}
\label{s:intro}

A Gabor system consists of functions of the form
\begin{equation}
\gnm(t) = e^{2\pi i mbt} g(t-na),\,\,\, n,m \in \Zst,
\label{gaborsystem}
\end{equation}
for some window $g \in \LtR$ and $a,b \in \Rst$.
The parameters $a$ and $b$ are the time- and frequency translation
parameters.

We say that the triple $(g,a,b)$ generates a {\em Gabor frame} for $\LtR$ 
for given $a,b$ if there exist constants ({\em frame bounds}) $A,B>0$ such
that
\begin{equation}
A \|f\|^2 \le \sum_{n,m \in \Zst} |\langle f, \gnm \rangle|^2 \le B\|f\|^2,
\label{gaborframe}
\end{equation}
for any $f \in \LtR$.

The Gabor frame operator $S$ is given by 
\begin{equation}
Sf = \sum_{n,m \in \Zst} \langle f, \gnm \rangle \gnm , \qquad f \in \LtR,
\label{frameop}
\end{equation}
and satisfies
\begin{equation}
\label{framebound}
 A I \le S \le B I, 
\end{equation}
where $I$ is the identity operator of $\LtR$.

If $(\gnm)_{n,m \in \Zst}$ is a Gabor frame for $\LtR$ then any $f$ in $\LtR$
can be written in the form
\begin{equation}
\label{framerep}
f = \sum_{n,m \in \Zst} \langle f,\gdnm \rangle \gnm =
\sum_{n,m \in \Zst} \langle f,\gnm \rangle \gdnm,
\end{equation}
where the {\em canonical dual frame} $(\gdnm)_{n,m \in \Zst}$ is given by
\begin{equation}
\label{dual}
\gdnm = e^{2\pi i mbt} \gd(t-na),\,\,\, n,m \in \Zst
\end{equation}
with $\gd = S^{-1} g$.

If $A=B$ the frame is called {\em tight}. In this case $S=\frac{1}{A}I$,
whence $\gd = \frac{1}{A}g$. We say that a frame is {\em normalized tight}
if in addition $A=B=1$.
Tight Gabor frames play an important role in signal processing and
communications. They appear for instance in the construction of paraunitary
modulated filter banks~\cite{CV98,BH98} as well as in the
construction of orthogonal frequency division multiplex (OFDM)
systems~\cite{Bol99,Str00} in wireless communications. Tight Gabor frames 
are also useful tools for the analysis of pseudodifferential 
operators~\cite{RT98}.

Given a frame $\gabframe$ with frame
operator $S$ a standard technique to construct a tight frame 
is the following. Set 
\begin{equation}
\label{tightframe}
\ho=\SQI g,
\end{equation}
then it is easy to see that $\tightframe$ is a
tight frame for $\LtR$. We will refer to $\ho$ as {\em canonical tight
window} and to $\tightframe$ as {\em canonical tight frame}.

It is known that the canonical dual $\gd$ minimizes
\begin{equation}
\label{min1}
\Big\|\frac{g}{\|g\|} - \frac{\gamma}{\|\gamma\|}\Big\|
\end{equation}
among all dual windows $\gamma$, cf.~\cite{Jan95}.
This minimality property of $\gd$ can be expressed more precisely
in the time-domain as: for a.e.~$t$ and all $\alpha>0, \beta>0$
(possibly depending on $t$) we have that
\begin{equation}
\label{min1a}
\sum_{n=-\infty}^{\infty} |\alpha g(t-na)-\beta \gamma(t-na)|^2
\end{equation}
is minimal among all dual $\gamma$ for $\gamma=\gd$.
Such minimality results also hold in the frequency domain and in the
Zak transform domain.

We will show in Section~\ref{s:tight} that a result similar to~\eqref{min1}
holds for tight Gabor frames.
More precisely, we will demonstrate that $\ho$ solves
\begin{equation}
\underset{h \in \LtR}{\argmin} \|g - h\|
\label{min}
\end{equation}
among all $h$ such that $(h,a,b)$ is a normalized tight frame.

This result is also important from a practical viewpoint, where one often 
selects a ``nice'' $g$ and constructs the tight frame by 
using~\eqref{tightframe} with the goal in mind that the resulting tight 
frame window should have similar properties as $g$.
Thus it is gratifying to know that $\ho$ is closest to $g$ (in the 
$\Ltsp$-sense) among all possible normalized tight frames $(h,a,b)$.
For conditions under which $\ho$ inherits the decay properties
of $g$ see~\cite{Wal92,Str98a,BJ00,Str00}.

It is sufficient to restrict our attention to normalized tight frames due
to following reason. Assuming that~\eqref{min} holds the reader may easily
convince herself/himself that the window $\tilde{h}^0$ that minimizes 
$\|g-h\|$ among all $h$ such that $(h,a,b)$ is a tight frame is given by
\begin{equation}
\label{min2}
\tilde{h}^0 = \frac{\langle \ho, g \rangle}{\langle \ho,\ho \rangle}\ho .
\end{equation}


For the special case that $(g,a=1,b=1)$ constitutes a frame it is not 
difficult to see that~\eqref{min} holds. In this case we can consider
the problem in the Zak transform domain with $\lambda=p=q=1$
(see Subsection~\ref{ss:rep}~\eqref{i14b}--\eqref{zak} for the 
definition of the Zak transform and the Zibulski-Zeevi matrices that are 
implicitly used in the next steps). Thus we have for any $f$ such that 
$(f,1,1)$ has a finite frame upper bound that
\begin{equation}
(Z_1 Sf)(t,\nu) = |(Z_1g)(t,\nu)|^2 (Z_1f)(t,\nu), \quad \aetn.
\label{b1}
\end{equation}
By functional calculus (cf.~Subsection~\ref{ss:rep}) we then have that
\begin{gather}
(Z_1 \ho)(t,\nu) = (Z_1 \SQI g)(t,\nu) =
\Big(|(Z_1 g)(t,\nu)|^2\Big)^{-\frac{1}{2}} (Z_1 g)(t,\nu) \\
=\frac{(Z_1 g)(t,\nu)}{|(Z_1 g)(t,\nu)|}, \quad \aetn.
\label{b2}
\end{gather}
Moreover, normalized tightness of $(h,1,1)$ is equivalently expressed as 
\begin{equation}
\label{b2a}
|(Z_1 h)(t,\nu)|=1,\quad  \aetn. 
\end{equation}
Now when we have a complex number $z$ then
\begin{equation}
\label{b3}
\underset{|\omega|=1}{\min} |z-\omega| = \Big|z - \frac{z}{|z|}\Big| = 
\Big| |z|-1\Big|.
\end{equation}
Thus we have that for a.e.~$t,\nu$ and all normalized tight $(h,1,1)$ 
\begin{equation}
\Big|(Z_1 g)(t,\nu) - (Z_1 \ho)(t,\nu)\Big|^2 \le 
\Big|(Z_1 g)(t,\nu) - (Z_1 h)(t,\nu)\Big|^2. 
\label{b4}
\end{equation}
Integrating the latter inequality over $[0,1)^2$ we get using unitarity
of $Z_1: \LtR \rightarrow \Ltsp([0,1)^2)$ that $\|g-\ho\|^2\le \|g-h\|^2$.

\bigskip

The rest of the paper is organized as follows.
Subsection~\ref{ss:rep} contains a brief review of various representations
of the Gabor frame operator in different domains. We will make use
of these representations in Section~\ref{s:tight} when we prove the 
minimization result~\eqref{min}. We also discuss some refinements and
extensions of this result. In Section~\ref{s:so} we present a
Wiener-Levy type theorem for Gabor frames. Finally, in Section~\ref{s:newton} 
we describe a Newton-type method for computing $\ho$ and
analyze the convergence properties of this method.\footnote{We were kindly 
informed by Hans G.~Feichtinger that he discovered the algorithm of 
Section~\ref{s:newton} independently of us, as early as 1995 
and used it in numerical experiments to compute
$\ho$ with satisfactory results.}

It is obvious that some or many of our results admit
generalization to more general lattices in $\Rst^2$ as sets of time-frequency
points to which the windows are shifted and to higher dimensions, but
we have not pursued this point here.

\subsection{Representations of the Gabor frame operator} \label{ss:rep}

We briefly review the representations of the Gabor frame operator in
the time domain, in the frequency domain, in the time-frequency domain,
and, for rational oversampling where $ab \in \Qst$, in the Zak transform 
domain. We also present a lemma concerning the approximation of the inverse 
square root of an infinite-dimensional matrix by finite-dimensional inverse
square roots.

Assume that $(g,a,b)$ has a finite frame upper bound $B$, and define
the Ron-Shen matrices
\begin{gather}
\label{i10}
\Mg(t):= \Big( g(t-na-k/b)\Big)_{k \in \Zst, n \in \Zst} ,\quad \aet,\\
\label{i11}
H_g(\nu):= \Big( \hat{g}(\nu-mb-l/a)\Big)_{l \in \Zst, m \in \Zst} ,\quad \aen.
\end{gather}
Now there holds that $\Mg(t)$ and $H_g(\nu)$ are for a.e.~$t$ and a.e.~$\nu$
bounded linear operators of $\ltZ$ (with operator norms 
$\le (bB)^{\frac{1}{2}}$ and $\le (aB)^{\frac{1}{2}}$). Now when
$(f,a,b)$ also has a finite frame upper bound, then there holds
\begin{gather}
\label{i12}
M_{Sf}(t) = \frac{1}{b} \Mg(t) \Mga(t) M_f(t) ,\quad \aet ,\\
\label{i13}
H_{Sf}(\nu) = \frac{1}{a} \Hg(\nu) \Hga(\nu) H_f(\nu) ,\quad \aen .
\end{gather}
For the case~\eqref{i11},~\eqref{i13} this can be found 
in~\cite{Jan98}, proof of Theorem~1.2.6 and Proposition~1.2.3, in the more
general context of shift-invariant systems, and, in somewhat different
form, in~\cite{RS95a} or in~\cite{RS97a}.

Next, still assuming that $(g,a,b)$ has a finite frame upper bound $B$,
we define the analysis operator with respect to the adjoint lattice
$\{(k/b, l/a) | k,l \in \Zst \}$ by
\begin{equation}
m \in \LtR \rightarrow \Ug m=\Big(\langle m,\gkl \rangle\Big)_{k,l\in \Zst}.
\label{i14}
\end{equation}
Now there holds that $\Ug$ is a bounded linear operator (operator norm
$\le(ab B)^{\frac{1}{2}}$). When $(f,a,b)$ also has a finite upper frame bound,
then there holds
\begin{equation}
U_{Sf} = \frac{1}{ab} \Ug \Uga U_f .
\label{i14a}
\end{equation}
These results can be found in~\cite{Jan95}, Propositions~2.7,~3.1, and
also in~\cite{RS97a} and~\cite{DLL95} (in somewhat different form).
Expression~\eqref{i14a} is sometimes called the Janssen representation of the 
frame operator.

Finally still assuming that $(g,a,b)$ has a finite frame upper bound $B$,
but now with $ab=p/q \in \Qst$ and $1 \le p \le q$, $\gcd(p,q)=1$, we
define the Zibulski-Zeevi matrices as
\begin{gather}
\Phi^g(t,\nu)=
p^{-\frac{1}{2}}\Big((Z_{\frac{1}{b}} g)(t-l\frac{p}{q},\nu+\frac{k}{p})\Big)
_{k=0,\dots,p-1;l=0,\dots,q-1} \label{i14b} \\
\Psi^g(t,\nu)=
p^{-\frac{1}{2}}\Big((Z_{a} g)(t-k\frac{q}{p},\nu-\frac{l}{q}) \Big)
_{k=0,\dots,p-1;l=0,\dots,q-1} \label{i14c} 
\end{gather}
for a.e.~$t,\nu \in \Rst$. Here we define for $\lambda >0$ the 
Zak transform $Z_{\lambda} f$ of an $f \in \LtR$ as
\begin{equation}
(Z_{\lambda}f)(t,\nu) = \lambda^{\frac{1}{2}}
\sum_{k=-\infty}^{\infty} f(\lambda (t-k)) e^{2\pi i k \nu},
\qquad \aetn \in \Rst.
\label{zak}
\end{equation}
Now there holds that $\Phi^g (t,\nu)$ and $\Psi^g (t,\nu)$ are for
a.e.~$t,\nu$ bounded linear mappings of $\Cst^p$ into $\Cst^q$
(matrix norms $\le B^{\frac{1}{2}}$). When $(f,a,b)$ also has a finite
frame upper bound, then there holds  
\begin{gather}
\Phi^{Sf}(t,\nu)= \Phi^g (t,\nu) \Big(\Phi^g (t,\nu)\Big)^{\ast}
\Phi^f (t,\nu), \qquad \aetn, \label{i15}  \\
\Psi^{Sf}(t,\nu)= \Psi^g (t,\nu) \Big(\Psi^g (t,\nu)\Big)^{\ast}
\Psi^f (t,\nu), \qquad \aetn. \label{i16}
\end{gather}
These results can be found in~\cite{Jan95b} and in~\cite{ZZ93}.

The frame operator $S$ is represented through the operators
\begin{gather}
S(t) := \frac{1}{b} \Mg(t)\Mga(t), \quad \hat{S}(\nu) := \frac{1}{a} \Hg(\nu)
\Hga(\nu),
\label{i17} \\
\underline{S} := \frac{1}{ab} \Ug \Uga, \label{i18}
\end{gather}
and
\begin{equation}
S_1(t,\nu) := \Phi^g(t,\nu) \Big(\Phi^g (t,\nu)\Big)^{\ast} ,\,\,\,
S_2(t,\nu) := \Psi^g(t,\nu) \Big(\Psi^g (t,\nu)\Big)^{\ast} 
\label{i19}
\end{equation}
in the following sense. When $A>0, B < \infty, g \in \LtR$, then
\begin{align}
& \text{$(g,a,b)$ is a Gabor frame with frame bounds $A, B$} \label{i20}\\
& \Longleftrightarrow \,\, A I \le S(t) \le B I , \aet \notag \\
& \Longleftrightarrow \,\, A I \le \hat{S}(\nu) \le B I , \aen \notag \\
& \Longleftrightarrow \,\, A I \le \underline{S} \le B I \notag \\
& \Longleftrightarrow \,\, A I \le S_1(t,\nu) \le B I , \aetn \notag \\
& \Longleftrightarrow \,\, A I \le S_2(t,\nu) \le B I , \aetn \notag .
\end{align}
Here the $I$'s denote identity operators of $\ltZ, \ltZ, 
\ltsp(\Zst \times \Zst), \Cst^{p}$, and $\Cst^p$, respectively. See the
references just given. 

An important point for the developments in this paper is that the relations
\eqref{i12}, \eqref{i13}, \eqref{i14a}, \eqref{i15}, \eqref{i16} can
be extended as follows. Assume that $(g,a,b)$ is a frame with frame 
bounds $A>0, B < \infty$, and let $\phi$ be a function analytic in an open
neighbourhood of $[A,B]$. When $(f,a,b)$ has a finite frame upper bound,
there holds
\begin{align}
M_{\phi(S)f}(t)&=\phi(S(t)) M_f(t), \qquad \aet , \label{i21} \\
H_{\phi(S)f}(\nu)&=\phi(\hat{S}(\nu)) H_f(\nu), \qquad \aen , \label{i22} \\
U_{\phi(S)f}&=\phi(\underline{S}) U_f, \label{i23} \\
\Phi^{\phi(S)f}(t,\nu)&=\phi(S_1(t,\nu))\Phi^f(t,\nu),\quad \aetn \label{i24}\\
\Psi^{\phi(S)f}(t,\nu)&=\phi(S_2(t,\nu)) \Psi^f(t,\nu),\quad \aetn. \label{i25}
\end{align}
Indeed, the aforementioned relations can be shown to hold by iteration when
$\phi$ is a polynomial, and then a standard argument in functional calculus
noting~\eqref{i20} then yields the result for general $\phi$ analytic
in an open neighbourhood of $[A,B]$. In the important special case that
$\phi(s) = s^{\alpha}$, one can do this extension explicitly by using
\begin{equation}
s^{\alpha} = \Big( \frac{B+A}{2} (1-v)\Big)^{\alpha} = 
\Big(\frac{B+A}{2}\Big)^{\alpha} \sum_{n=0}^{\infty} \binom{\alpha}{n} (-v)^n ,
\label{i26}
\end{equation}
where $v= 1-2s(B+A)^{-1}$ has modulus $\le (B-A)(B+A)^{-1}$ when
$s \in [A,B]$.

\bigskip
A further result that we shall use is that for any frame $(g,a,b)$ and any
normalized tight frame $(h,a,b)$ we have
\begin{equation}
\label{i27}
\langle g, \SI g \rangle = \|h\|^2 = ab.
\end{equation}
See~\cite{DLL95,Jan95,RS97a}.

In the sequel we will also make use of the following result.
\begin{lemma}
\label{le:kant}
Let $T=(T_{kl})_{k,l \in \Zst}$ be a hermitian positive definite biinfinite
matrix with $c_1 I \le T \le c_2 I$ . Define $T_N = (T_{kl})_{|k|,|l| \le N}.$
Then $\Tni \rightarrow \Ti$ strongly for $N \rightarrow \infty$.
\end{lemma}

\begin{proof}
Without loss of generality assume $c_2 <1$, otherwise we can
always set $T' = \frac{1}{c_2}T$. Note that
\begin{equation}
\label{dummy2}
\Ti = \sum_{k=0}^{\infty} c_k (I-T)^k,
\end{equation}
where
$c_k=\binom{-\frac{1}{2}}{k}$ and
$\Tni = \sum_{k=0}^{\infty} c_k (I_N-T_N)^k $.
Define
$\TKi = \sum_{k=0}^{K} c_k (I-T)^k$
and
$\TKni = \sum_{k=0}^{K} c_k (I_N-T_N)^k$.
Observe that 
\begin{equation}
c_1 I \le \|T_N \|\le c_2 I \quad \text{uniformly in $N$}
\label{kant1}
\end{equation}
and
\begin{equation}
T_N \rightarrow T \quad \text{strongly for $N \rightarrow \infty$.}
\label{kant2}
\end{equation}
For fixed $N$ there holds $\TKni \rightarrow \Tni$ as $K \rightarrow \infty$.
Furthermore, for each $K$ we have that
$\TKni \rightarrow \TKi$ for $N \rightarrow \infty$. This together with
properties~\eqref{kant1} and~\eqref{kant2} implies that 
$\Tni \rightarrow \Ti$ strongly for $N \rightarrow \infty$.
\end{proof}

\section{Canonical tight Gabor frames and the \\
Fan-Hoffman inequality} 
\label{s:tight}

\begin{theorem}
\label{th:fan}
Let $(g,a,b)$ be a frame and let $\ho=\SQI g$ be the canonically 
associated tight frame generating window. Then for any $h \in \LtR$ for 
which $(h,a,b)$ is a normalized tight Gabor frame there holds
\begin{equation}
\|g - \ho \| \le \|g - h\| \le \|g + \ho\|.
\label{fanhoffman}
\end{equation}
\end{theorem}

In Subsections~\ref{ss:ronshen},~\ref{ss:janssen}, and~\ref{ss:zz}
we will present different proofs of this inequality, making use
of different representations of the Gabor frame operator.

Inequality~\eqref{fanhoffman} can be considered as a Fan-Hoffman type 
inequality~\cite{FH55}. The Fan-Hoffman inequality can be formulated as
follows. Let the matrix $M \in \Cst^{m \times n}, m \ge n$ have the polar 
decomposition $M=UH$, where $U^{\ast} U=I_n$ and $H$ is hermitian 
positive semidefinite. Then for any $m \times n$ matrix $W$ with
$W^{\ast} W=I_n$ 
\begin{equation}
\|M-U\| \le \|M - W\| \le \|M+U\|
\label{fh1}
\end{equation}
for any unitarily invariant norm. See~\cite{FH55} for a proof for
$M \in \Cst^{n \times n}$ and~\cite{Hig86} for a discussion of
the case $M \in \Cst^{m \times n}$.

An immediate consequence of the Fan-Hoffman inequality is
the following result.
Consider a frame for $\Cst^n$ consisting of $m$ vectors
$\phi_k, k=1,\dots,m$ with frame operator $S$. Then the tight frame 
$(\psi_k^0)_{k=1}^m$ that minimizes
\begin{equation}
\sum_{k=1}^{m} \|\phi_k - \psi_k\|^2
\label{sum}
\end{equation}
among all tight frames $(\psi_k)_{k=1}^m$ in $\Cst^n$ is given by
\begin{equation}
\label{dummy1}
\psi_k^0 = \frac{\trace(S^{\frac{1}{2}})}{n} \SQI \phi_k .
\end{equation}
This can be easily seen 
by considering the vectors $\phi_k$ as rows of an $m \times n$ matrix $M$
and using~\eqref{fh1} with the Frobenius norm, see also Section~1.5 
in~\cite{Str99}.

The case where $n=m$ and $(\phi_k)_{k=1}^n$ is an orthonormal basis has 
been investigated by L\"owdin in the context of quantum chemistry~\cite{GL91}. 
The solution -- known under the name {\em L\"owdin orthogonalization} --
is of course a special case of~\eqref{dummy1}.

Both, the Fan-Hoffman inequality and the L\"owdin orthogonalization, 
have been extended to infinite dimensions, see e.g.~Chapter~VI
in~\cite{GK69}, or~\cite{AEG80,FPT00}. These extensions usually require that
the involved operators satisfy some Hilbert-Schmidt-type condition
in order that
\begin{equation}
\sum_{k \in \Zst} \|\phi_k - \psi_k\|^2
\label{low1}
\end{equation}
and the operator version of~\eqref{fh1} are finite.
A different notion of closeness of frames, based on
work of Paley and Wiener~\cite{PW87}, has been used in~\cite{Bal99}.

For Gabor frames the expression $\sum_{m,n \in \Zst} \|\gnm - \hnm\|^2$
is of course not finite for $g \neq h$. But since $\|\gnm - \hnm\|=\|g-h\|$ 
for all $n,m$ we can restrict ourselves to the minimization of $\|g-h\|$.
The minimization result for tight Gabor frames presented in this
paper can also be considered directly in the context of L\"owdin 
orthogonalization, due to the following reason. The results in
Subsection~\ref{ss:rep} imply that if $(h,a,b)$ is a tight frame for
$\LtR$ then $(h,1/b,1/a)$ is an orthogonal basis for its closed
linear span. Thus the statement that $\ho$ minimizes $\|g-h\|$ among
all normalized tight frames $(h,a,b)$ is equivalent to saying
that $\ho$ minimizes $\|g-h\|$ among all $h$ such that
$(h,1/b,1/a)$ is an orthonormal basis for the closed linear span
of $(\gnm)_{n,m}$.

\subsection{The problem in the time domain and in the frequency domain} 
\label{ss:ronshen}

\subsubsection{The problem in the time domain} \label{sss:time}

For the first proof of Theorem~\ref{th:fan} we consider the Gabor frame 
operator in the time domain using Ron-Shen matrices.

We consider matrices 
\begin{equation}
\Mg(t) = \Big(g(t-na-\kb)\Big)_{k \in \Zst, n \in \Zst} \qquad \aet,
\label{rs1}
\end{equation}
with row index $k$ and column index $n$. We have by functional calculus,
see~\eqref{i21}
\begin{equation}
\Mho(t) = \Big(\frac{1}{b} \Mg(t) \Mga(t)\Big)^{-\frac{1}{2}} \Mg(t)\qquad \aet.
\label{rs2}
\end{equation}
The condition of normalized tightness can be expressed equivalently as,
see~\eqref{i20}
\begin{equation}
\Mh(t)\Mha(t)=b I \qquad \aet. 
\label{rs3}
\end{equation}
Note also that by normalized tightness $\|h\|^2=ab$ (see~\eqref{i27}). Thus
\begin{equation}
\|h-g\|^2 = ab + \|g\|^2 - 2 \Re \langle h, g \rangle,
\label{rs4}
\end{equation}
and
\begin{equation}
\|h+g\|^2 = ab + \|g\|^2 + 2 \Re \langle h, g \rangle,
\label{rs4a}
\end{equation}
and so we are to maximize $\Re \langle h,g \rangle$ over all $h$ 
satisfying~\eqref{rs3}.

Let $K=1,2,\dots$ and consider the sections
\begin{equation}
\Mgk(t) = \Big(g(t-na-\kb)\Big)_{|k|\le K, n\in \Zst}.
\label{rs5}
\end{equation}
$\Mgk$ consists of $(2K+1)$ rows and a biinfinite number of columns.
Note that
\begin{equation}
\Mgk(t)\Mgka(t)=\Big(\big(\Mg(t)\Mga(t)\big)_{k,l}\Big)_{|k|\le K,|l|\le K}.
\label{rs6}
\end{equation}
Thus we have that for any allowed $h$ there holds
\begin{equation}
\Mgk(t)\Mgka(t)=b\, I_{2K+1}.
\label{rs7}
\end{equation}
We are going to maximize, for fixed $t \in \Rst$, the quantity
\begin{equation}
\trace\big(\Mk \Mgka(t)\big),
\label{rs8}
\end{equation}
where $\Mk=(M_{kn})_{|k|\le K,n \in \Zst}$ is a matrix satisfying
\begin{equation}
\Mk \Mka = b\, I_{2K+1}.
\label{rs9}
\end{equation}
Noting that the matrix in~\eqref{rs6} is positive definite (as a
$(2K+1)\times (2K+1)$ section of the positive definite operator
$\Mg(t) \Mga(t)$ since $(g,a,b)$ is a frame for $\LtR$) we do a
singular value decomposition of $\Mgk(t)$ according to
\begin{equation}
\Mgk(t) = \sum_{j=0}^{2K} \sigma_j u_j v_j^H,
\label{rs10}
\end{equation}
where
\begin{gather}
v_j = (v_{jn})^T_{n \in \Zst} \quad \in \ltZ \notag  \\
u_j = [u_{j,-K},\dots, u_{j,K}]^T \quad \in \Cst^{2K+1} \notag
\end{gather}
are orthonormal, $\sigma_j^2$ is the $j$-th eigenvalue of $\Mgk(t)\Mgka(t)$
and $u_j$ is the corresponding eigenvector. Thus
\allowdisplaybreaks
\begin{equation}
\big(\Mgk(t)\Mgka(t) \big) u_j = \sigma_j^2 u_j
\label{rs11}
\end{equation}
and
\begin{equation}
v_j = \frac{1}{\sigma_j} \Mgka (t) u_j.
\label{rs12}
\end{equation}
Now we compute
\begin{gather}
|\trace \big(\Mk \Mgka(t)\big)|=
\big|\trace \big(\Mk \sum_{j=0}^{2K} \sigma_j v_j u_j^H\big)\big|\notag \\
=\big|\trace \big(\sum_{j=0}^{2K} \sigma_j M_k v_j u_j^H\big)\big|
=\big|\trace \big(\sum_{j=0}^{2K} \sigma_j u_j^H M_k v_j\big) \big| \notag \\
=\big|\sum_{j=0}^{2K} \sigma_j \langle \Mk v_j, u_j \rangle \big|
\le \sum_{j=0}^{2K} \sigma_j \|\Mk v_j\| \| u_j \| \notag \\
\le b^{\frac{1}{2}} \sum_{j=0}^{2K} \sigma_j =
 b^{\frac{1}{2}} \trace \Big(\Mgk(t) \Mgka(t) \Big)^{\frac{1}{2}}.
\end{gather}
Taking 
\begin{equation}
\Mk(t) = \Big(\frac{1}{b} \Mgk(t)\Mgka(t)\Big)^{-\frac{1}{2}} \Mgk(t)
\label{rs13}
\end{equation}
we also see that the maximum value
\begin{equation}
b^{\frac{1}{2}} \trace \Big(\Mgk(t)\Mgka(t)\Big)^{\frac{1}{2}}
\label{rs14}
\end{equation}
is assumed by $\Mk(t)$ in~\eqref{rs13}.

Next we observe that for an allowed $h$ we have
\begin{equation}
\trace [\Mhk(t)\Mgka(t)] = \sum_{k=-K}^{K} \sum_{n=-\infty}^{\infty}
h(t-na-\kb)g^{\ast}(t-na-\kb).
\label{rs15}
\end{equation}
Taking the average over $k$ and integrating over $t \in [0,a)$, we thus see
that
\begin{equation}
\begin{array}{c}
\langle h,g \rangle = \frac{1}{2K+1}\int \limits_{0}^{a}\sum_{k=-K}^{K} 
\sum_{n=-\infty}^{\infty}h(t-na-\kb)g^{\ast}(t-na-\kb)\, dt \\
= \frac{1}{2K+1}\int \limits_{0}^{a} 
\trace [\Mhk(t)\Mgka(t)]\, dt .
\end{array}
\label{rs16}
\end{equation}
The right-hand side of~\eqref{rs16} has modulus $\le$ what one would obtain
by replacing $\Mhk(t)$ by $\Mk(t)$ in~\eqref{rs16}. 
Due to Lemma~\ref{le:kant} we have that
\begin{equation}
\underset{K \toinf} \lim \Mk(t) = 
\underset{K \toinf}\lim \Big(\frac{1}{b}\Mgk(t)\Mgka(t)\Big)^{-\frac{1}{2}}
\Mgk(t) = M_{\ho}(t)
\label{rs17}
\end{equation}
in the strong operator topology. This implies that
\begin{equation}
\underset{K \toinf} \lim \frac{1}{2K+1} \int \limits_{0}^{a} 
\trace [\Mk(t)\Mgka(t)] dt = \langle \ho, g \rangle,
\label{rs18}
\end{equation}
and the proof of inequality~\eqref{fanhoffman} is complete.

\bigskip
Actually we have shown more than~\eqref{fanhoffman}. Consider the
$a$-periodic $\Lsp^{\infty}$-functions
\begin{equation}
\label{rs19}
f(t)=\sum_{n=-\infty}^{\infty}h(t-na)g^{\ast}(t-na)=\langle h,g \rangle_a (t),
\end{equation}
where we have used the bracket product notation 
$\langle \cdot,\cdot\rangle_a$ as in~\cite{CL99}. 
That $f \in \Lsp^{\infty}$ follows from the frame bound conditions,
together with the considerations in Subsection~\ref{ss:rep}. Consider the case 
that $ab=p/q$ with integer $p,q, 1 \le p\le q, \gcd(p,q)=1$. Then for 
a.e.~$t$ there holds
(see the right-hand side of~\eqref{rs15})
\allowdisplaybreaks
\begin{gather}
\label{rs20}
\frac{1}{2K+1} \sum_{k=-K}^{K} f(t-k/b) \\
= \frac{1}{2K+1} \sum_{k=-K}^{K} f(t-\frac{q}{p} a) \notag \\
\underset{[k=lp+r]}{=}
\frac{1}{2K+1}\sum_{r=0}^{p-1}\sum_{l;-K\le lp+r\le K}f(t-lqa-r\frac{q}{p}a)
\notag \\
\underset{\text{[using the $a$-periodicity of $f$}]}{=} 
\frac{1}{2K+1}\sum_{r=0}^{p-1} \#_{l;-K\le lp+r\le K}\cdot f(t-r\frac{q}{p}a)
\notag \\
\approx \frac{1}{2K+1} \cdot \frac{2K+1}{p} \sum_{r=0}^{p-1}f(t-r/b)
\notag \\
=\frac{1}{p} \sum_{r=0}^{p-1}\langle h,g \rangle_a (t-r/b). \notag
\end{gather}
Thus there holds, more precisely, that for a.e.~$t$
\begin{equation}
\label{rs21}
\frac{1}{p} \sum_{r=0}^{p-1} \langle h,g \rangle_a (t-r/b)
\end{equation}
is maximized among all normalized tight $(h,a,b)$ by $\ho$. And 
then, of course, integrating over $t \in [0,1/b)$ gives the result.

In case that $ab \notin \Qst$, it can be shown that
\begin{equation}
\frac{1}{2K+1} \sum_{k=-K}^{K} f(\cdot - k/b) \rightarrow \langle h,g \rangle
\label{rs22}
\end{equation}
in the mean square sense when $K \toinf$ (here $f$ is as in~\eqref{rs19}).
Hence in this case we do not get a more precise result of the type given
above for rational $ab$.

\subsubsection{The problem in the frequency domain} \label{sss:frequency}

We now present a proof of Theorem~\ref{th:fan} where we consider the
Gabor frame operator in the frequency domain, using frequency domain
Ron-Shen matrices. We thus consider
\begin{equation}
\Hg(\nu) = \Big(\hat{g}(\nu -mb-l/a)\Big)_{l\in \Zst, m \in\Zst} ,\quad
\aen ,
\label{f1}
\end{equation}
and we have by functional calculus
\begin{equation}
H_{\ho}(\nu) = \Big(\frac{1}{a} \Hg(\nu) \Hga(\nu)\Big)^{-\frac{1}{2}}
\Hg(\nu), \quad \aen,
\label{f2}
\end{equation}
while the condition of normalized tightness can be expressed equivalently as
\begin{equation}
H_h(\nu) H_h^{\ast}(\nu) = a I , \quad \aen.
\label{f3}
\end{equation}
We therefore find ourselves in a situation that is similar to what we had
in Subsection~\ref{sss:time} in all respects, except that the Ron-Shen 
matrices involve now the Fourier transform of the windows, rather than the
windows themselves. Accordingly, we know that $\|\hat{h}-\hat{g}\|^2$
is minimal among all $h$ such that $(h,a,b)$ is normalized tight for $h=\ho$.
By Parseval's theorem we have that $\|\hat{h}-\hat{g}\|^2=
\|h-g\|^2$, and the proof is complete. We also get a frequency domain
sharpening of the result for $ab\in \Qst$ in a similar manner as the one
obtained at the end of Subsection~\ref{sss:time}. This latter sharpening
reads that for a.e.~$\nu$
\begin{equation}
\label{f4}
\frac{1}{p} \sum_{s=0}^{p-1} \langle \hat{h},\hat{g} \rangle_b (\nu-s/a)
\end{equation}
is maximal among all $h$ such that $(h,a,b)$ is normalized tight for $h=\ho$.

\bigskip

The frequency domain result can be generalized to some extent to 
shift-invariant systems
\begin{equation}
(\gsnm)_{n,m\in\Zst} = \Big(g_m(\cdot - na)\Big)_{n,m\in\Zst},
\label{f5}
\end{equation}
where $g_m \in \LtR$ and $a>0$. Following the presentation of the
Ron-Shen theory~\cite{RS95a} given in~\cite{Jan98}, Subsection~1.2,
we now consider the frame operator
\begin{equation}
Sf = \sum_{n,m} \langle f , \gsnm \rangle \gsnm , \qquad f \in \LtR,
\label{f6}
\end{equation}
and the Ron-Shen matrices
\begin{equation}
H_g(\nu)=\Big(g_m(\nu-l/a)\Big)_{l\in \Zst,m\in\Zst}, \quad \aen .
\label{f7}
\end{equation}
There holds for $A>0, B < \infty$ that
\begin{equation}
AI \le S \le B I \Leftrightarrow A I \le \frac{1}{a} H_g(\nu)
H_g^{\ast}(\nu) \le B I , \enspace \aen.
\label{f8}
\end{equation}
The canonical dual system $(\gd_{nm})_{n,m \in \Zst}$ and the
canonically associated tight system $(\hosnm)_{n,m \in \Zst}$ are given by
\begin{equation}
\Big(\gd_m(\cdot-na)\Big)_{n,m\in\Zst}, \enspace
\Big(\homm(\cdot-na)\Big)_{n,m\in\Zst}
\label{f9}
\end{equation}
with $\gd_m = S^{-1} g_m, \homm=\SQI g_m$. Here we note that $S$
in~\eqref{f6} commutes with all shift operators $f \rightarrow f(\cdot-na),
n \in \Zst$. The conditions of duality of $(\gamma_{nm})_{n,m \in \Zst}$
and normalized tightness of $(h_{nm})_{n,m \in \Zst}$ can be equivalently 
expressed as
\begin{equation}
\label{f10}
H_g(\nu)H_{\gamma}^{\ast}(\nu) = aI ,\quad  H_h(\nu) H_h^{\ast} (\nu) = aI.
\end{equation}
We have a functional calculus in terms of the matrices
$a^{-1}H_g(\nu)H_g^{\ast}(\nu)$ by the formula
\begin{equation}
H_{\phi(S)g}(\nu) = \phi\Big(\frac{1}{a}H_g(\nu)\Hga(\nu)\Big)\Hg(\nu),\quad \aen,
\label{f11}
\end{equation}
where $\phi$ is a function analytic in an open neighbourhood of $[A,B]$.
We have set here
\begin{equation}
\phi(S)g = \Big((\phi(S)g_m)(\cdot-na)\Big)_{n,m\in \Zst},\quad \aen.
\label{f12}
\end{equation}
In particular there holds
\begin{equation}
H_{\ho}(\nu) = \Big(\frac{1}{a}\Hg(\nu)\Hga(\nu)\Big)^{-\frac{1}{2}}
\Hg(\nu), \enspace \aen.
\label{f13}
\end{equation}
Now the developments as given in Subsection~\ref{sss:time} can be mimicked
to a large extent. Considering truncated matrices $H_{g,K}$ as 
in~\eqref{rs5}, and denoting the Frobenius norm by $\| . \|_F$,
the minimization of (for a.e.~$\nu$)
\begin{equation}
\|H_K - H_{g,K}(\nu)\|^2_F = \sum_{|l|\le K}\sum_{m=-\infty}^{\infty}
|h_{lm}-\hat{g}(\nu-l/a)|^2
\label{f14}
\end{equation}
over all $(2K+1) \times \infty$ matrices $H_K=(h_{lm})_{|l|\le K,m\in\Zst}$
satisfying
\begin{equation}
H_K H_K^{\ast} = a I_{2K+1}
\label{f15}
\end{equation}
yields that the minimizing $H_K$ is given by
\begin{equation}
H_K = \Big(\frac{1}{a} H_{g,K}(\nu)H_{g,K}^{\ast}(\nu)\Big)^{-\frac{1}{2}}
H_{g,K}(\nu).
\label{f16}
\end{equation}
And, again by Lemma~\ref{le:kant}, this $H_K \rightarrow H_{\ho}$ as
$K \toinf$ in the strong operator topology.

\subsection{The problem in the time-frequency domain} \label{ss:janssen}

We next present a proof of Theorem~\ref{th:fan} where we consider the
Gabor frame operator in the time-frequency domain, using analysis operators
$\Ug$ as defined in~\eqref{i14}
based on the adjoint lattice $\{(k/b, l/a) | k,l \in \Zst\}$.

Note that $\Ug \Uga$ has the matrix
representation
\begin{equation}
\Ug\Uga=\Big(\langle \gklp,\gkl \rangle \Big)_{k,l \in \Zst; k',l' \in\Zst}.
\label{j3}
\end{equation}

Normalized tightness is equivalently expressed as
\begin{equation}
\Uh \Uha = ab \, I .
\label{j5}
\end{equation}
The canonical tight frame generating $\ho$ is given by
\begin{equation}
U_{\ho} = \Big( \frac{1}{ab} \Ug \Uga\Big)^{-\frac{1}{2}} \Ug
\label{j6}
\end{equation}
by functional calculus. We also note that
\begin{equation}
\Uh \Uga=\Big(\langle \gklp, \hkl \rangle \Big)_{k,l\in \Zst;k',l'\in\Zst}.
\label{j7}
\end{equation}
Now consider the sections
\begin{equation}
\Ugkl : m \in \LtR \rightarrow 
\Big( \langle m, \gkl \rangle \Big)_{|k|\le K, |l| \le L} 
\in \Cst^{2K+1} \times \Cst^{2L+1}
\label{j8}
\end{equation}
with adjoints 
\begin{equation}
\Ugkla : \underline{c} \in \Cst^{2K+1}\times \Cst^{2L+1} \rightarrow 
\sum_{\stackrel{|k|\le K}{|l| \le L}} c_{kl} \gkl
\label{j9}
\end{equation}
and 
\begin{equation}
\Ugkl \Ugkla = \Big( (\Ugkl \Ugkla)_{k,l;k'l'}\Big)_{|k|\le K,|l|\le L;
|k'|\le K, |l'|\le L}.
\label{j10}
\end{equation}
This we also do with allowed $h$'s, see~\eqref{j5}, and we thus see that
\begin{equation}
\Uhkl \Uhkla = ab \, I_{(2K+1)(2L+1)}.
\label{j11}
\end{equation}
Now consider the maximization of 
\begin{equation}
\Big| \trace \Big[U_{K,L} \Ugkla \Big] \Big|
\label{j12}
\end{equation}
over all 
\begin{equation}
U_{K,L} = \Big( U_{kl;k'l'}\Big)_{|k|\le K,|l|\le L; |k'|\le K, |l'|\le L}
\label{j13}
\end{equation}
such that $U_{K,L} U_{K,L}^{\ast} = ab I_{(2K+1)(2L+1)}$. This yields
the maximum value
\begin{equation}
\trace \Big[ (ab)^{\frac{1}{2}} \Big(\Ugkl \Ugkla \Big)^{\frac{1}{2}}\Big],
\label{j14}
\end{equation}
assumed by
\begin{equation}
U_{K,L} = \Big(\frac{1}{ab} \Ugkl \Ugkla \Big)^{-\frac{1}{2}} \Ugkl .
\label{j15}
\end{equation}
Note that from~\eqref{j10}
\begin{equation}
\langle g,h \rangle = \frac{1}{(2K+1)(2L+1)} \trace \Big[\Uhkl \Ugkla \Big].
\label{j16}
\end{equation}
By Lemma~\ref{le:kant} $U_{K,L}$ in~\eqref{j15} converges strongly to
$U_{\ho,K,L}$, hence
\begin{equation}
\langle g,\ho \rangle = \underset{K,L \toinf}{\lim}
\frac{1}{(2K+1)(2L+1)} \trace \Big[U_{K,L} \Ugkla \Big].
\label{j17}
\end{equation}
Therefore, $|\langle g,h \rangle| \le \langle g,\ho \rangle$.

\subsection{The problem in the Zak transform domain} \label{ss:zz}

We proceed by presenting a proof of Theorem~\ref{th:fan} by considering
the Gabor frame operator in the Zak transform domain, using Zibulski-Zeevi
matrices. The proof in this section is the simplest among all presented
proofs, however it only applies for rational oversampling.

Assume $ab=p/q, 1 \le p \le q, \gcd(p,q)=1$. By~\cite{Jan98}, 
Sec.~1.5,~(1.5.11) we have
\begin{align}
\langle h,g \rangle & =  p \int \limits_{0}^{q^{-1}} \int \limits_{0}^{p^{-1}}  
\sum_{k=0}^{p-1} \sum_{l=0}^{q-1} \Phi_{kl}^h (t,\nu) 
\big(\Phi^g_{kl} (t,\nu)\big)^{\ast} dt\, d\nu \notag \\
& =  p \int \limits_{0}^{q^{-1}} \int \limits_{0}^{p^{-1}}  
\trace \big[ \Phih \Phiga \big] dt\, d\nu \notag .
\end{align}
Here $\Phi^f(t,\nu)$ is the Zibulski-Zeevi matrix as defined
in~\eqref{i14b}.
The condition of normalized tightness is expressed equivalently as
\begin{equation}
\Phih \Phiha = I_p \qquad \aetn.
\label{zz4}
\end{equation}
The maximum of 
\begin{equation}
\Big| \trace \big[\Phi \Phiga \big] \Big|
\label{zz5}
\end{equation}
over all $p \times q$ matrices $\Phi$ with $\Phi \Phi^{\ast} = I_p$ is
equal to 
\begin{equation}
\trace \Big[\Big(\Phig \Phiga \Big)^{\frac{1}{2}}\Big]
\label{zz6}
\end{equation}
and assumed by
\begin{equation}
\Phi (t,\nu) : =\Big(\Phig \Phiga \Big)^{-\frac{1}{2}} \Phig .
\label{zz7}
\end{equation}
The right-hand side of~\eqref{zz7} equals $\Phi_{\ho} (t,\nu)$
by functional calculus for Zibulski-Zeevi matrices.

Hence we see that $\ho$ indeed maximizes $\langle h,g \rangle$,
and that it does so in a pointwise manner in terms of Zibulski-Zeevi
matrices.

The proof using the matrices $\Psi^g$ as defined in~\eqref{i14c} is almost
identical to the one for the matrices $\Phi^g$ and is therefore
left to the reader.

\section{A Wiener-Levy theorem for Gabor frames} \label{s:so}

In~\cite{Wal92,Str98a,BJ00,Str00} it has been investigated under which
conditions the canonical tight window inherits the decay properties
of the window $g$. Loosely spoken $\ho$ inherits the decay properties of
$g$ whenever this is true for $\gd$.
We give a precise mathematical formulation of this observation.
We start with the following result from the theory of Banach 
algebras~\cite{Gar66}.
\begin{proposition}
\label{prop:gardner}
Let $\Ban$ be a Banach algebra, and $x$ an element of $\Ban$ with
real-valued spectrum $\sigma (x)>0$. Then there exists in $\Ban$ a unique
square root of $x$ with real-valued positive spectrum.
\end{proposition}

This implies the following: Assume that the properties of $g$ imply that
$S$ belongs to some Banach algebra $\Ban$. Whenever $S^{-1} \in \Ban$,
it follows that $\SQI \in \Ban$. 

Actually for rational values of $ab$ a more general result holds.
\begin{theorem}
\label{th:wienerlevy}
Assume that $ab=p/q, 1\le p,q \in \Zst, \gcd(p,q)=1$. Assume that
$(g,a,b)$ is a frame, where $g$ satisfies Tolimieri and Orr's 
condition A (cf.~\cite{TO95})
\begin{equation}
\sum_{k,l} |\langle g,\gkl \rangle| < \infty.
\label{tol1}
\end{equation}
Let $\Win$ be a time-frequency shift-invariant Banach space of windows
$w$ such that $(w,a,b)$ has a finite frame upper bound when $w \in \Win$.
Finally, let $\phi$ be analytic in an open neighbourhood of $[A,B]$,
where $A>0, B<\infty$ are the frame bounds of $(g,a,b)$.  
Then $\phi(S)w \in \Win$ for any $w \in \Win$.
\end{theorem}

For the proof of Theorem~\ref{th:wienerlevy} we first need a matrix-version
of the classical Wiener-Levy theorem. We introduce the following
\begin{definition}
A set $U$ of functions $f(t)=\sum_{k\in \Zst} a_k(f) e^{2\pi i kt}$ 
is said to have an absolutely uniformly convergent
Fourier series if for every $\eps >0$ there is a $K >0$ such that
\begin{equation}
\sum_{|k| \ge K} |a_k(f)| \le \eps, \qquad \text{for all $f \in U$}.
\label{uniform}
\end{equation}
\end{definition}


\begin{proposition}
Let $H(t)$ be a periodic hermitian matrix-valued function with absolutely
convergent Fourier series, i.e.,  
$H(t)=\sum_{k\in \Zst} H_k e^{2\pi i kt}, H_k=((H_k)_{m,n})_{m,n=1}^{N}$
and 
\begin{equation}
\label{fnorm}
\|H\|_{\cal F}:= \sum_{m,n=1}^{N}\sum_{k\in \Zst}|(H_k)_{m,n}| < \infty.
\end{equation}
Let $A, B \in \Rst, A < B$ be such that all eigenvalues of all $H(t)$
are in $[A,B]$. If the function $\phi(z)$ is analytic in an open neighbourhood
of $[A,B]$ then $\phi(H(t))$ also has an absolutely
convergent Fourier series.
\label{prop:levy}
\end{proposition}
\begin{proof}
We have the Dunford representation~\cite{DS63}
\begin{equation}
\phi(H(t)) = \frac{1}{2\pi i} \oint \limits_{\Gamma} \phi(z)
\big(zI - H(t)\big)^{-1}\, dz,
\label{dunford}
\end{equation}
where $\Gamma$ is a closed contour containing $[A,B]$ in its interior.
We write
\begin{equation}
\big(zI-H(t)\big)^{-1}=\det \big(zI-H(t)\big)^{-1}\adj\big(zI-H(t)\big).
\label{det}
\end{equation}
The entries of $zI-H(t)$ have Fourier series absolutely uniformly
convergent in $z \in \Gamma$. Hence it follows easily that 
$\adj\big(zI-H(t)\big)$ and $\det \big(zI-H(t)\big)$, as sums and products
of the entries of $zI-H(t)$, have Fourier series that are absolutely 
uniformly convergent in $z \in \Gamma$. Furthermore, there are constants
$C_1 >0, C_2 < \infty$ such that
\begin{equation}
C_1 \le \big| \det\big(zI-H(t)\big)\big| \le C_2,\quad z \in \Gamma, t\in \Rst.
\label{bound}
\end{equation}
Thus by the uniform Wiener $1/f$-theorem~\cite{CCJ00}
we have that $\det \big(zI-H(t)\big)^{-1}$ has a Fourier series that is
absolutely uniformly convergent in $z \in \Gamma$. This implies that
$\|(zI-H)^{-1}\|_{\cal F}$ is uniformly bounded in $z \in \Gamma$, by $C$, 
say. Hence
\begin{equation}
\|\phi(H)\|_{\cal F} \le \frac{1}{2\pi} C 
\oint \limits_{\Gamma} | \phi(z)| |dz| < \infty,
\label{phibound}
\end{equation}
where the boundedness of the integral in~\eqref{phibound} follows
from the continuity of $\phi$. The proof is complete.
\end{proof}

The knowledgeable reader will have no difficulties to extend the argument 
concerning the uniform Wiener $1/f$-theorem given in~\cite{CCJ00} to the case 
that we are considering functions of two variables with absolutely uniformly 
convergent Fourier series. With the two-dimensional version of the uniform 
Wiener $1/f$-theorem at hand it is then easy to extend 
Proposition~\ref{prop:levy} to bivariate hermitian 
matrix-valued functions $H(t,\nu)$.

\noindent {\bf Proof of Theorem~\ref{th:wienerlevy}: }
Let $w \in \Win$. By functional calculus we have 
\begin{equation}
\phi(S)w=\sum_{k,l}\Big(\phi(\frac{1}{ab} \Ug \Uga)\Big)^{\ast}_{0,0;k,l}
w_{k/b,l/a}, \label{tol2}
\end{equation}
with convergence of the right-hand series at least in weak $\Ltsp$-sense.
It is therefore sufficient to show that
\begin{equation}
\sum_{k,l}\Big|\Big(\phi(\frac{1}{ab} \Ug \Uga)\Big)^{\ast}_{0,0;k,l}\Big|
< \infty
\label{tol3}
\end{equation}
since $w_{k/b,l/a} \in \Win$ for all $k,l \in \Zst$.

To show~\eqref{tol3} we take $w=g$ in~\eqref{tol2}. By biorthogonality we
have
\begin{equation}
\Big(\phi(\frac{1}{ab} \Ug \Uga)\Big)^{\ast}_{0,0;k,l}=
\frac{1}{ab} \langle \phi(S)g, \gamma^{0}_{k/b,l/a} \rangle.
\label{tol4}
\end{equation}
From Prop.~1.1 in~\cite{Jan95b} we have for any $f, h \in \LtR$
\begin{gather}
\label{tol5}
\sum_{k,l}|\langle f,h_{k/b,l/a} \rangle | < \infty \Longleftrightarrow \\
\Phi^f(t,\nu) \big(\Phi^h(t,\nu)\big)^{\ast}\,\,
\text{has an abs.~conv.~Fourier series.} \notag
\end{gather}
Now take 
\begin{equation}
f=\phi(S)g, \qquad h=\gd = S^{-1}g.
\label{tol6}
\end{equation}
By functional calculus in the Zak transform domain, we have
\begin{gather}
\label{tol7}
\Phi^f(t,\nu)\big(\Phi^h(t,\nu)\big)^{\ast}= \\
=\phi\Big[\Phi^g(t,\nu)\big(\Phi^g(t,\nu)\big)^{\ast}\Big] \Phi^g(t,\nu) 
\Big(\big[\Phi^g(t,\nu)\big(\Phi^g(t,\nu)\big)^{\ast}\Big]^{-1} 
\Phi^g(t,\nu) \Big)^{\ast} \notag \\
=\phi \Big[\Phi^g(t,\nu)\big(\Phi^g(t,\nu)\big)^{\ast}\Big].\notag
\end{gather}
By the bivariate version of Proposition~\ref{prop:levy} we have that
the right-hand side of~\eqref{tol7} has an absolutely convergent
Fourier series. Then by~\eqref{tol4}, \eqref{tol5}, \eqref{tol6} we
see that~\eqref{tol3} holds, and the proof is complete. \QED

\medskip
\noindent 
{\bf Example:}
Assume that $ab =p/q, 1\le p,q \in \Zst, \gcd(p,q)=1$. If $g$ is in the 
Feichtinger-algebra $\SOsp$~\cite{FG96}, then by Theorem~\ref{th:wienerlevy}
$\ho \in \SOsp$.

\section{A Newton method for the computation of tight Gabor frames} 
\label{s:newton}

For applications and design purposes it is useful to have a fast algorithm
for computing the canonical tight window $\ho$. It is clear that in the case 
of integer oversampling this can be done via the Zak transform. However for
non-integer oversampling the frame operator can at best only be
block-diagonalized. This means that at some point we would have to 
compute the inverse square root of a matrix, which is computationally 
much more expensive than the inversion of a matrix. In this section we present 
a Newton-type algorithm that allows a very efficient calculation of $\ho$.

Newton's method is well-known in the context of computing matrix-sign
functions and square roots of finite-dimensional matrices~\cite{Hig86a,KL92}.
Let $M$ be an $n \times n$ matrix with positive real-valued eigenvalues
and polar decomposition $M=UH$.
Set $M_0=M, X_0=M$, and define the iterations
\begin{align}
&M_{k+1}=\frac{1}{2} (M_k + [M_k^{\ast}]^{-1}M), \qquad k=0,1,\dots
\label{newton1}\\
&X_{k+1}=\frac{1}{2} (X_k + [X_k^{\ast}]^{-1}), \qquad k=0,1,\dots.
\label{newton2}
\end{align}
Then $M_k$ converges to $M^{\frac{1}{2}}$ and $X_k$ converges to $U$,
see~\cite{Hig86a,KL92}. 

In Section~\ref{s:tight} we have seen that $\ho$ minimizes $\|g-h\|$ among
all windows $h$ for which $(h,a,b)$ is a normalized tight frame. Clearly in the 
same way $\ho$ minimizes $\|\gd-h\|$. Also recall that $\gd$ minimizes
\begin{equation}
\Big\| \frac{g}{\|g\|} - \frac{\gamma}{\|\gamma\|}\Big\|
\label{minimumgd}
\end{equation}
among all dual windows $\gamma$. This suggests to compute $\ho$
iteratively by setting $g_0=g$ and defining the 
iteration
\begin{equation}
g_{k+1} = \frac{1}{2}\Big(\frac{g_k}{\|g_k\|} + \frac{\gd_k}{\|\gd_k\|}\Big), 
\qquad k=0,1,\dots,
\label{newton3}
\end{equation}
where $\gd_k$ is the dual window associated with the frame $(g_k,a,b)$.
Accordingly, we denote the frame bounds of the system $(g_k,a,b)$ by $A_k$ 
and $B_k$. Iteration~\eqref{newton2} can be interpreted as a scaled Newton 
iteration. The proof of convergence of Newton's method as described 
in~\eqref{newton1} and~\eqref{newton2} relies on finite-dimensional methods 
(see~\cite{Hig86}), which cannot be extended straightforwardly to infinite 
matrices or operators.

\subsection{Convergence of the algorithm} \label{ss:conv}

In this section we analyze the convergence behavior of the
iteration~\eqref{newton3}. We need the following 
\begin{lemma}
\label{le:analytic}
Assume that $(g,a,b)$ is a frame, and let $\phi$ be a function analytic in
an open neighbourhood of $[A,B]$, where $A$ and $B$ are the optimal frame
bounds for $(g,a,b)$, and assume that $\phi(s)>0, s \in [A,B]$. Then
$(\phi(S)g,a,b)$ is a frame with frame bounds
\begin{equation}
\label{le1}
\underset{A\le s\le B}{\min} s \phi^2(s) , \enspace 
\underset{A\le s\le B}{\max} s \phi^2(s) .
\end{equation}
Furthermore $(g,a,b)$ and $(\phi(S)g,a,b)$ have the same canonically
associated tight frame generating window, viz.~$\ho$.
\end{lemma}
\begin{proof}
We have for any $f \in \LtR$ that
\begin{equation}
\label{le2}
\sum_{n,m} \langle f, (\phi(S)g)_{nm} \rangle (\phi(S)g)_{nm} =
\phi(S) [\sum_{n,m} \langle (\phi^{\ast}(S)f,g_{nm} \rangle]
= \phi(S) S \phi^{\ast}(S) f .
\end{equation}
Here we have used that $S$, and hence $\phi(S)$, commutes with all
relevant time-frequency shifts. Thus the frame operator corresponding to
$(\phi(S)g,a,b)$ is given by $\phi(S) S \phi^{\ast}(S)$. By the spectral
mapping theorem we have that
\begin{equation}
\sigma [\phi(S) S \phi^{\ast}(S)]= \{\phi(s)s\phi^{\ast}(s)=s|\phi(s)|^2
\,\,\big| s \in \sigma(S)\},
\label{le3}
\end{equation}
and this gives the first part of the result. 

Next we compute the tight frame generating window $h^{\phi}$ canonically
associated to $(\phi(S)g,a,b)$ according to
\begin{equation}
\label{le4} 
h^{\phi}=\Big(\phi(S) S \phi(S)\Big)^{-\frac{1}{2}} \phi(S)g 
=\SQI (\phi^2(S))^{-\frac{1}{2}} \phi(S)g = \SQI g =\ho.
\end{equation}
Here we have used that $\phi(s)>0, s \in [A,B]$ so that $\phi^{\ast}(S)
=\phi(S)$ and $\phi(S) > 0$ so that
$\Big((\phi(S))^2\Big)^{-\frac{1}{2}}=\phi(S)$.
\end{proof}

\begin{theorem}
\label{th:conv}
Let $(g,a,b)$ be a frame for $\LtR$ with optimal frame
bounds $A,B$ and canonical tight window $\ho$. Set $g_0=g$ and define
\begin{equation}
\label{iter}
g_{k} = \frac{1}{2}\Big(\frac{g_{k-1}}{\|g_{k-1}\|} + 
\frac{\gd_{k-1}}{\|\gd_{k-1}\|}\Big), 
\qquad k=1,2,\dots.
\end{equation}
then $g_{k}$ converges quadratically to $\frac{1}{\sqrt{ab}} \ho$.
\end{theorem}

\proof
Set
\begin{equation}
g_1=\phi(S) g \,:\, \phi(x) = \alpha + \beta x^{-1},
\label{n00}
\end{equation}
where for simplicity we have denoted $\alpha_k=\frac{1}{\|g_k\|}, 
\beta_k=\frac{1}{\|\gd_k\|}$ and $\alpha=\alpha_0, \beta=\beta_0$. Note that
\begin{equation}
A \le \frac{\beta}{\alpha} = \frac{\|g\|}{\|S^{-1}g\|} \le B.
\label{n0}
\end{equation}
Using Lemma~\ref{le:analytic} we obtain that
\begin{gather}
A_1=\underset{x \in [A,B]}{\min} (\alpha^2 x + 2\alpha\beta +\beta^2 x^{-1}),
\label{n3}\\
B_1=\underset{x \in [A,B]}{\max} (\alpha^2 x + 2\alpha\beta +\beta^2 x^{-1}).
\label{n3a}
\end{gather}
As to $A_1$, we note that $\alpha^2 x +2\alpha \beta +\beta^2 x^{-1}$
assumes its minimum value $4\alpha \beta$ at $x=\beta/\alpha \in [A,B]$,
see~\eqref{n0}. Hence $A_1 = 4\alpha \beta$, and we furthermore note that
\begin{equation}
4\alpha \beta = \frac{1}{\|g\| \|S^{-1}g\|} \ge 
\frac{1}{|\langle g, S^{-1} g \rangle |} = \frac{1}{ab},
\label{n3b}
\end{equation}
see~\eqref{i27}. Hence $A_1 \ge \frac{1}{ab}$, no matter what $g$ is.
As to $B_1$ we note that $\alpha^2 x + 2\alpha\beta +\beta^2 x^{-1}$ is 
strictly convex on $[A,B]$ and assumes its minimum at $x=\beta/\alpha$.
Hence 
\begin{equation}
B_1=\underset{x=A \,\text{or}\,B}{\max}(\alpha^2 x+2\alpha\beta+\beta^2 x^{-1}).
\label{n4}
\end{equation}
Therefore we find that
\begin{align}
\frac{A_1}{B_1} & = \min \Big\{ 
\frac{4 \alpha \beta}{\alpha^2 A+2\alpha \beta+\beta^2 A^{-1}},
\frac{4 \alpha \beta}{\alpha^2 B+2\alpha \beta+\beta^2 B^{-1}}\Big\}
\notag \\ 
& = \min \Big\{ \frac{4A \frac{\beta}{\alpha}}{(A+\frac{\beta}{\alpha})^2},
\frac{4B \frac{\beta}{\alpha}}{(B+\frac{\beta}{\alpha})^2}\Big\}.
\label{n5}
\end{align}
\begin{figure}
\begin{center}
\epsfig{file=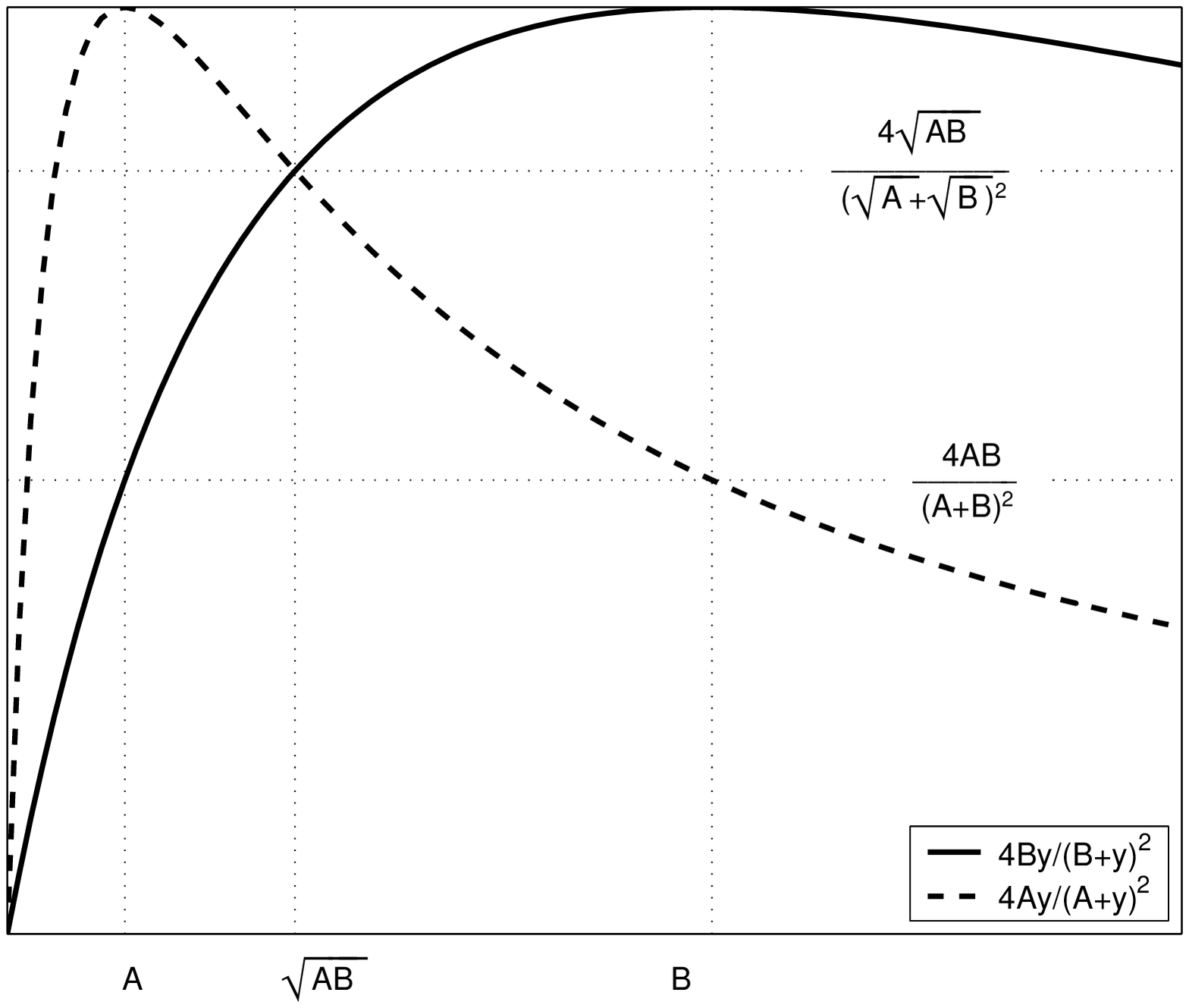,width=100mm,height=80mm}
\caption{}
\label{fig:opt}
\end{center}
\end{figure}
The graph in Figure~\ref{fig:opt} shows the two functions
$f_1(y)=\frac{4Ay}{(A+y)^2}$, $f_2(y)=\frac{4By}{(B+y)^2}$ for
$y=\frac{\beta}{\alpha}$. It is obvious that
\begin{equation}
\frac{4 AB}{(A+B)^2} \le
\frac{A_1}{B_1} \le
\frac{4 \sqrt{AB}}{(\sqrt{A}+\sqrt{B})^2} \, .
\label{n6}
\end{equation}
The above argument applies to any $g_{k-1}$ as in~\eqref{iter}. In
particular, we have for $k=1,2,\dots$
\begin{equation}
1 \ge \frac{A_k}{B_k} \ge \frac{4 A_{k-1} B_{k-1}}{(A_{k-1}+B_{k-1})^2}
= \frac{4 \frac{A_{k-1}}{B_{k-1}}}{(1+\frac{A_{k-1}}{B_{k-1}})^2}
\label{n7}
\end{equation}
Consider the recursion
\begin{equation}
C_k = \frac{4C_{k-1}}{(1+C_{k-1})^2}, \qquad k=1,2,\dots
\label{n8}
\end{equation}
with $C_0 \in (0,1)$. The recursion~\eqref{n8} converges {\em quadratically}
to 1, since the function
\begin{equation}
f(y)=\frac{4y}{(1+y)^2},\qquad y \in (0,1]
\label{n9}
\end{equation}
increases on $(0,1]$, satisfies $f(0)=0, f(1)=1$, $f(y)>y, y\in (0,1)$,
and $f'(1)=0, f''(1)=-\frac{1}{2}$.
Evidently when we take $C_0=A/B$, we obtain that 
\begin{equation}
1 \ge \frac{A_k}{B_k} \ge C_k ,\qquad k=1,2,\dots
\label{n10}
\end{equation}
and $A_k/B_k \rightarrow 1$, quadratically. 

It follows from~\eqref{iter} and~\eqref{n00} that the $g_k$ are all of the 
form $\phi_k (S) g$. Hence by Lemma~\ref{le:analytic} they all have
$\ho$ as their canonically associated tight frame generator. At the
same time, denoting by $S_k$ the frame operator of the system
$(g_k,a,b)$, we have that this canonically associated tight frame
generator is given by $\SQI_k g_k$. Hence
 \begin{equation}
\ho = S^{-\frac{1}{2}}_k g_k ,\qquad k=1,2,\dots .
\label{n11}
\end{equation}
Since $A_k/B_k \rightarrow 1$, we easily see from
\begin{equation}
g_k=\frac{1}{2}\Big(\frac{g_{k-1}}{\|g_{k-1}\|}+
    \frac{S^{-1}_{k-1} g_{k-1}}{\|S^{-1}_{k-1} g_{k-1}\|} \Big)
\label{n12}
\end{equation}
that $\|g_k\| \rightarrow 1$ as $k \toinf$. Also $\|\ho\|=\sqrt{ab}$
by~\eqref{i27}, and it then follows from~\eqref{n11} and the
fact that $A_k/B_k \rightarrow 1$ that
\begin{equation}
A_k, B_k \rightarrow \frac{1}{ab}, \qquad k \toinf .
\label{n13}
\end{equation}
All these convergences are quadratic.
Consequently $S^{-\frac{1}{2}}_k \rightarrow \sqrt{ab}\, I$ quadratically, 
which completes the proof.  \QED

\smallskip

Considering~\eqref{n5} and the left-hand-side of~\eqref{n6} it is not
difficult to check that quadratic convergence of the algorithm holds for 
other choices of the scaling parameters. Writing the iteration as
\begin{equation}
g_{k+1} = \frac{1}{2}\Big( \alpha_k g_k + \beta_k \gd_k \Big), \qquad
k=0,1,\dots
\label{n13a}
\end{equation}
it is easy to see that the algorithm converges quadratically 
to a scaled version of $\ho$ if $\frac{\beta_k}{\alpha_k} \in [A,B]$.
The right-hand side in~\eqref{n6} shows that the convergence cannot be
expected to be faster than quadratic. Indeed, although the function 
\begin{equation}
\label{dumm5}
k(y)=\frac{4y^{\frac{1}{2}}}{(1+y^{\frac{1}{2}})^2}=f(y^{\frac{1}{2}}),
\end{equation}
(see~\eqref{n9}) is considerably flatter than $f(y)$ at
$y=1$, we have that $k(1)=1, k'(1)=0, k''(1)=-\frac{1}{8}$, hence we have
at best that $1-\frac{A_k}{B_k}=\frac{1}{16}(1-\frac{A_{k-1}}{B_{k-1}})^2$
compared to $1- \frac{A_k}{B_k} = \frac{1}{4}(1-\frac{A_{k-1}}{B_{k-1}})^2$
for $f(y)$ in~\eqref{n9}.
This improves convergence behavior with respect to the constant, but not
with respect to the order. However, since from a numerical viewpoint
improving the constant is of considerable importance, we want to 
find out which choice of $\alpha_k, \beta_k$ leads to the optimal
constant.

Optimizing~\eqref{n5} with respect to $\frac{\beta}{\alpha}$
yields that the optimum choice for the scaling parameters is
\begin{equation}
\alpha_k = \beta_k^{-1} = (A_k B_k )^{-\frac{1}{4}},
\label{n14}
\end{equation}
see also the graph in Figure~\ref{fig:opt}.
The same optimal scaling parameter arises in the
finite-dimensional Newton method, cf.~\cite{Hig86} as well
in Balan's construction of ``nearest tight frames'', see
Section~3 in~\cite{Bal99}.
It is not difficult to see that for optimal scaling we get
\begin{equation}
\label{n14a}
A_1=4, \qquad B_1=\frac{(\sqrt{A}+\sqrt{B})^2}{\sqrt{AB}},
\end{equation}
hence $\frac{A_1}{B_1}$ takes on the maximum value in the right-hand side
of~\eqref{n6}.

The choice~\eqref{n14} however is practically not feasible, since it requires
the computation (or estimation) of the frame bounds $A_k, B_k$ in
each iteration step, which is in general computationally expensive.

An alternative is to find the scaling parameter $\alpha$  that minimizes
\begin{equation}
\label{dummy6}
\|\alpha g - \alpha^{-1} \gd\|.
\end{equation}
A simple calculation shows that 
in this case 
\begin{equation}
\alpha=\frac{\|\gd\|}{\|g\|},
\label{frob}
\end{equation}
which is similar to Frobenius-norm scaling in~\cite{KL92} (hence we will 
refer to it as Frobenius-norm scaling) 
and closely related to the scaling used in Theorem~\ref{th:conv}.

Next we estimate how the frame bounds get ``tighter''
from iteration $k$ to $k+1$ when using the scaling proposed in
Theorem~\ref{th:conv}.

\begin{corollary}
\label{cor:est}
Under the assumptions of Theorem~\ref{th:conv} the 
frame bounds $A_k, B_k$ of $(g_k,a,b)$ can be estimated 
recursively by
\begin{equation}
A_k \ge \frac{1}{ab} \frac{2\sqrt{A_{k-1}B_{k-1}}}{A_{k-1}+B_{k-1}},\qquad
B_k \le \frac{1}{ab} \frac{(A_{k-1}+B_{k-1})^2}{4A_{k-1}B_{k-1}}.
\label{frameest}
\end{equation}
\end{corollary}

\begin{proof}
It is sufficient to consider the case 
$k=1$. We begin with the observation that for $S$ satisfying 
$A I \le S \le B I$ where $A>0, B < \infty$ and $f \in \LtR$ there holds
\begin{equation}
\frac{\langle S^{-1}f ,f \rangle}{\|f\| \|S^{-1}\|} 
\ge \frac{2\sqrt{AB}}{A+B}.
\label{n15}
\end{equation}
Inequality~\eqref{n15} is equivalent to the operator version of
the Kantorovich inequality~\cite{Kan52}.
For the reader's convenience we include
a short proof. Inequality~\eqref{n15} follows from
\begin{gather}
(A+B)^2 \langle \SI x,x \rangle^2 - 4AB \|x\|\|\SI x\|^2 = \notag \\
= \Big[ 2AB \|\SI x\|^2 - (A+B) \|\SQI x\| \Big]^2 +\notag  \\
+ 4AB \langle (BI-S)(I-A\SI) \SQI x, \SQI x \rangle \|\SI x\| \ge 0,
\end{gather}
since $ (BI-S)(I-A\SI)$ is a non-negative operator.

As a consequence we have now immediately that
\begin{equation}
A_1 = \frac{1}{\|g\|\|S^{-1}g\|} = \frac{1}{ab} 
\frac{\langle S^{-1}g ,g \rangle}{\|g\| \|S^{-1}g\|} \ge \frac{1}{ab}
\frac{2\sqrt{AB}}{A+B},
\label{n16}
\end{equation}
which is the first inequality in~\eqref{frameest}.

From~\eqref{n6} we have
\begin{equation}
B_1 \le \frac{(A+B)^2}{4AB} A_1 = \frac{1}{ab} \frac{(A+B)^2}{4AB}
\frac{\langle S^{-1}g ,g \rangle}{\|g\| \|S^{-1}g\|} 
\ge \frac{1}{ab} \frac{(A+B)^2}{4AB},
\label{n17}
\end{equation}
and this is the second inequality in~\eqref{frameest}.
\end{proof}

Both estimates in~\eqref{frameest} are sharp. This can be seen
by considering the following example. Assume $a=b=1$.
Let $\eps >0$ be arbitrarily small, and take $g$ such that
$|Z g| = B^{\frac{1}{2}}, 1/|Z g| = |ZS^{-1}g|=B^{-\frac{1}{2}}$ on a set of 
measure $1-\eps$ in the unit square while
$|Z g| = A^{\frac{1}{2}}, 1/|Z g| = |ZS^{-1}g|=A^{-\frac{1}{2}}$ on
the complement of this set. Then we have $\|g\|\approx B^{\frac{1}{2}}$,
$\|S^{-1} g\|\approx B^{-\frac{1}{2}}$, whence
\begin{gather}
B_1 = \underset{x = A,B}{\max} \Big(\frac{x}{4\|g\|^2} +
\frac{1}{2 \|g\| \|S^{-1}g\|} + \frac{x^{-1}}{4\|S^{-1}g\|^2} \Big)\\
\approx  \underset{x = A,B}{\max} \Big(\frac{x}{4B} +
\frac{1}{2} + \frac{x^{-1}}{4B^{-1}} \Big) = \frac{(A+B)^2}{4AB}.
\label{n18}
\end{gather}
Hence the second inequality in~\eqref{frameest} is sharp.

As to sharpness in the first inequality in~\eqref{frameest}, we consider
a $g$ that has a Zak transform with $Zg = B^{\frac{1}{2}}$ on a set
$\subset [0,1)^2$ of measure $t \in [0,1]$ and $Zg = A^{\frac{1}{2}}$
on the complementary set. Then
\begin{gather}
\label{n19}
A_1 = \frac{1}{\|g\| \|S^{-1}g\|} = \frac{1}{\|Zg\|\|Z S^{-1}g\|}\\
=\Big[(A+t(B-A))(\frac{1}{A}+t(\frac{1}{B}-\frac{1}{A}))\Big]^{-\frac{1}{2}}\\
=\Big(1+t \frac{(B-A)^2}{AB} - t^2 \frac{(B-A)^2}{AB}\Big)^{-\frac{1}{2}}
\end{gather}
Taking $t=\frac{1}{2}$ here, we have 
\begin{equation}
A_1 = \Big(1+\frac{(B-A)^2}{4AB}\Big)^{-\frac{1}{2}}=\frac{2\sqrt{AB}}{A+B},
\label{n20}
\end{equation}
which is the right-hand side of the first inequality in~\eqref{frameest}.

We compute $Z g_1$ for this last example. We find
\begin{equation}
Z g_1 = 
\begin{cases}
\frac{1}{2}\Big(\frac{B}{A+t(B-A)}\Big)^{\frac{1}{2}} +
\frac{1}{2}\Big(\frac{A}{B-t(B-A)}\Big)^{\frac{1}{2}}
& \text{when $Zg=B^{\frac{1}{2}}$}, \\
\frac{1}{2}\Big(\frac{A}{A+t(B-A)}\Big)^{\frac{1}{2}} +
\frac{1}{2}\Big(\frac{B}{B-t(B-A)}\Big)^{\frac{1}{2}}
& \text{when $Zg=A^{\frac{1}{2}}$}.
\end{cases}
\label{cases}
\end{equation}
The difference of upper and lower number in~\eqref{cases},
ignoring terms of order $(2t-1)^n (B-A)^m$ with $n+m> 4$
(there are no terms with $n+m=4$), is
\begin{equation}
\label{diff}
(1-2t)\frac{(\sqrt{B}-\sqrt{A})^2}{2\sqrt{AB}}.
\end{equation}
E.g., when $t=\frac{1}{2}$ the difference is zero.
This example also illustrates that the algorithm may converge even
much better than predicted by the theoretical estimates presented above.

\subsection{Numerical simulations} \label{ss:ex}

The numerical simulations presented in this section are representative
for a large number of situations. We consider a well-localized window
$g$ and oversampling rate $\frac{1}{ab} \in [1,2]$. Such a setup is typical
for oversampled modulated filter banks as well as for OFDM
systems (where $k$-times oversampling corresponds to $k$-times 
undersampling). We choose $g(t)=2^{\frac{1}{4}} e^{-\pi |t|^2}$ with 
$a=b=\sqrt{\frac{15}{16}}$ resulting in an oversampling rate of $16/15$. 
For the numerical implementation we follow the finite-dimensional
model in~\cite{Str99}, Sec.1.6.

\begin{figure}
\begin{center}
\epsfig{file=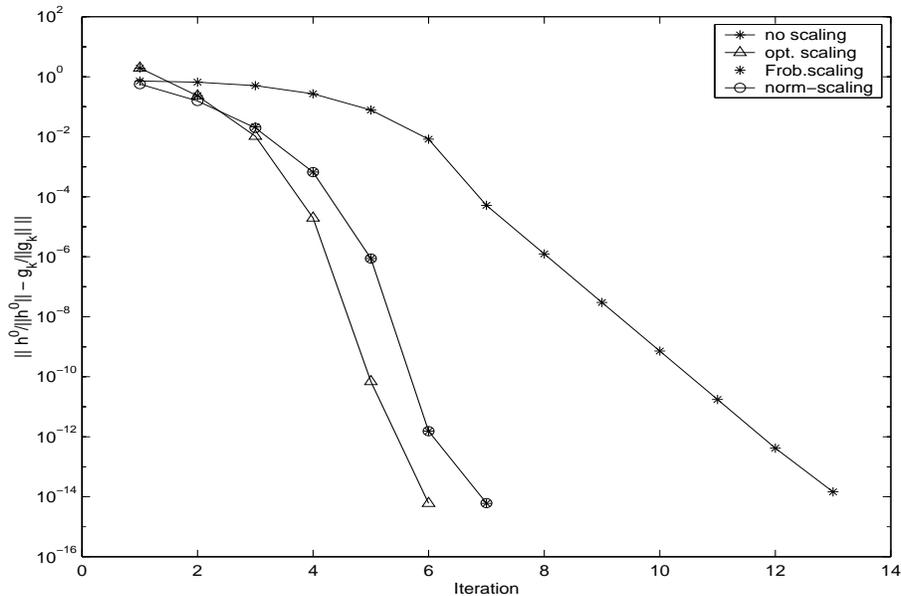,width=120mm,height=80mm}
\caption{Comparison of different scalings for Newton's method.}
\label{fig:ex1}
\end{center}
\end{figure}

In the first experiment we consider different types of scaling for
iteration~\eqref{n13a}. We compare the Newton method with 
{\em norm-scaling} as  used in Theorem~\ref{th:conv} to optimal scaling 
(cf.~\eqref{n14}), Frobenius-norm scaling
(cf.~\eqref{frob}) and to the unscaled Newton iteration
(i.e., $\alpha_k=\beta_k=1$).
We precompute $\ho$, calculate in each iteration the normalized
error $\|g_{k}/\|g_k\| - \ho/\|\ho\| \|$ and terminate the iterations
when the accuracy is within $10^{-14}$.
The results are shown in~Figure~\ref{fig:ex1}.

As expected, optimal scaling results in fastest convergence. Norm-scaling 
and Frobenius-norm scaling yield almost identical convergence, slightly
slower than optimal scaling. Recall however, that the costs
for computing the optimal scaling parameters are much larger than one
iteration step for the other two scaling schemes. Thus the faster convergence 
comes at a prohibitive high price. Newton's method without scaling requires 
significantly more iterations in order to achieve the same approximation error.
In terms of overall computational costs Newton's method with
norm-scaling or Frobenius-norm scaling are clearly the most
efficient methods.

An alternative to the algorithm presented in Theorem~\ref{th:conv} may
be to compute first iteratively $\SQI$ and at the end calculate
$\ho=\SQI g$. Efficient iterative methods to compute the (inverse) square root
of a positive definite matrix have been proposed in~\cite{She91,Lak95}.
These methods are different from Newton's method as described
in~\eqref{newton1}. 

\begin{figure}
\begin{center}
\epsfig{file=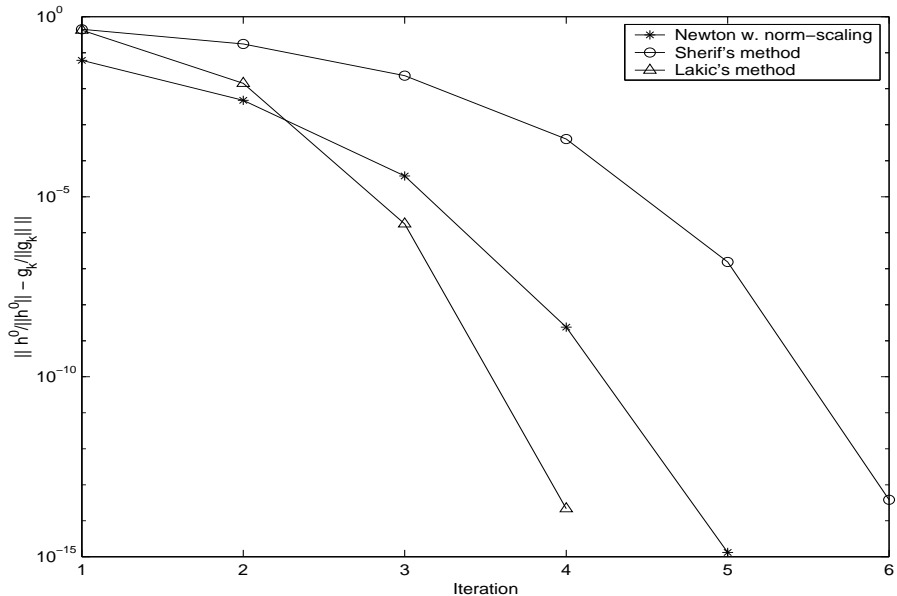,width=120mm,height=80mm}
\caption{Comparison of different methods to compute the tight window $\ho$.}
\label{fig:ex2}
\end{center}
\end{figure}

Let $M$ be a symmetric positive definite matrix. Then Sherif's method
to compute $M^{-\frac{1}{2}}$ is given by (cf.~\cite{She91}):
\begin{align}
& M_0  = I \notag \\
& M_{k+1} =  2 M_k (I + M M_k^2)^{-1} , \qquad k=0,1,\dots
\label{sherif}
\end{align}
This method yields quadratic convergence.
The method due to Lakic achieves even a cubic rate of convergence.
Its iteration rule is (cf.~\cite{Lak95}):
\begin{align}
& M_0  =I \notag \\
& M_{k+1}  = \frac{1}{3}M_k(I+8(I+3M M_k^2)^{-1}), k=0,1,\dots
\label{lakic}
\end{align}
Lakic's method yields usually slightly faster convergence than our method,
whereas Sherif's method is somewhat slower than ours, see
Figure~\ref{fig:ex2}. For an overall comparison of the efficiency of
iterative methods we have to take into account the computational
effort needed for one iteration step. The number of flops for the
matrix-multiplications and inversions in~\eqref{n12}, \eqref{sherif}, 
and~\eqref{lakic} depends on the chosen representation of the frame
operator. Note that for the proposed method this representation
has to be computed for each iteration (since we are dealing with a
new frame in each iteration), which is not the case for the
other two methods. However Newton's method needs only one matrix inversion per
iteration, whereas the other two methods require three matrix-matrix
multiplications in addition to one matrix inversion.

Fortunately the computation of the various representations of the
frame operator can in general be done very efficiently. For the
representation in the time-domain this amounts to simple data addressing,
in the other three domains FFT-based algorithm can be used
(possibly combined with data addressing).
Thus, as long as the computation of this representation is cheaper
than the inversion of the resulting matrix, the method 
proposed in Theorem~\ref{th:conv} outperforms
the other two methods in terms of computational efficiency.

\section*{Acknowledgement} 

The authors wish to thank Hans G.~Feichtinger for fruitful discussions
on the topic of this paper.

\end{document}